# ON UNIVERSAL ESTIMATES FOR BINARY RENEWAL PROCESSES

By Gusztáv Morvai [1] and Benjamin Weiss

*MTA-BME Stochastics Research Group and Hebrew University of Jerusalem*

A binary renewal process is a stochastic process $\{X_n\}$ taking values in $\{0,1\}$ where the lengths of the runs of 1's between successive zeros are independent. After observing $X_0, X_1, \ldots, X_n$ one would like to predict the future behavior, and the problem of universal estimators is to do so without any prior knowledge of the distribution. We prove a variety of results of this type, including universal estimates for the expected time to renewal as well as estimates for the conditional distribution of the time to renewal. Some of our results require a moment condition on the time to renewal and we show by an explicit construction how some moment condition is necessary.

## 1. Introduction.

The classical binary renewal process is a stochastic process $\{X_n\}$ taking values in $\{0,1\}$ where the lengths of the runs of 1's between successive zeros are independent. These arise, for example, in the study of Markov chains since the return times to a fixed state form such a renewal process; cf. [7]. (More details on this will be given in the next section.) In many applications, the occurrences of a zero, which represent the failure times of some system which is renewed after each failure, are of importance and so the problem arises of estimating when the next failure will occur; cf. Example 12.13 in [8].

Our purpose in this paper is to investigate the possibility of giving a universal estimator at time $n$ for the residual waiting time to the next zero in the binary renewal process $\{X_n\}$. Let $\{p_k\}_{k=0}^{\infty}$ be the conditional probability that a run of $k$ 1's follows a given 0 event. This distribution describes completely the renewal process as a two-sided stationary process. In order that the probability of $X_0 = 0$ be nonzero it is necessary that $\mu = \sum_{k=0}^{\infty} kp_k < \infty$

[1]Supported by the Bolyai János Research Scholarship during the second revision of this paper.

*AMS 2000 subject classifications.* 60G25, 60K05.
*Key words and phrases.* Prediction theory, renewal theory.

and then $P(X_0 = 0) = 1/(1+\mu)$ is positive. (This relation between the mean of the conditional renewal distribution and the stationary probability of the renewal event is well known in ergodic theory as Kac's formula for the expected return time to a set, and in probability theory cf. [7], Chapter XIII and [28], Section I.2.c.) If the process distribution is known, then after observing $X_0, X_1, \ldots, X_n$ one may give a consistent estimator for the expected value of residual waiting time to the occurrence of the next zero as

$$\mu_L = \frac{\sum_{k=L}^{\infty}(k-L)p_k}{\sum_{k=L}^{\infty} p_k}$$

if there is at least one zero among the values of $X_0, X_1, \ldots, X_n$ and the last zero occurs at moment $X_{n-L} = 0$. (Indeed, if $X_{n-L} = 0$ and for all $n - L < i \le n$, $X_i = 1$, then for $k \ge L$ the probability that for all $n+1 \le i < n - L + k + 1$, $X_i = 1$, and $X_{n+k-L+1} = 0$ is $\frac{p_k}{\sum_{i=L}^{\infty} p_i}$.) We denote this $L$ by $\tau(X_0, X_1, \ldots, X_n)$. Similarly we define $\tau = \tau(X_{-\infty}^0)$ as that $t \ge 0$ such that $X_{-t} = 0$ and $X_i = 1$ for all $-t < i \le 0$. It is clear from the stationarity that $P(\tau = L)$ is proportional to $\sum_{k=L}^{\infty} p_k$ and thus for the finiteness of the unconditional expectation of the residual waiting time we would have to demand that $\sum_{k=0}^{\infty} k^2 p_k < \infty$. We shall not assume this since we are interested primarily in the conditional expectations and with probability 1 for $n$ sufficiently large at least one of the $X_i = 0$, for $0 \le i \le n$ and for any fixed value of $\tau(X_0, \ldots, X_n) = L \le n$ the expected residual waiting time is $\mu_L < \infty$. This of course is well known in the classical analysis of renewal processes. In the spirit of recent investigations into universal estimators for various features of stationary processes (see [1, 2, 3, 6, 10, 13, 14, 17, 23, 24, 25, 29, 30]) we take up here the problem of how well can we do when all that we know is that the binary process $\{X_n\}$ is in fact a renewal process. The fact that we are trying to estimate the time to next occurrence of zero rather than $X_{n+1}$ takes us out of the framework of previous investigations. In earlier works such as [11] attention is restricted to those renewal processes which arise from Markov chains with a finite number of states. In that case the probabilities $p_k$ decay exponentially and one can use this information in trying to find not only the distribution but even the hidden Markov chain itself. Since we are considering the general case where the number of hidden states might be infinite, this exponential decay no longer holds in general and the problem becomes much more difficult.

For the estimator itself it is most natural to use the empirical distribution observed in the data segment $X_0, X_1, \ldots, X_n$. However, if there were an insufficient number of occurrences of 1-blocks of length at least $\tau(X_0, X_1, \ldots, X_n)$, then we do not expect the empirical distribution to be close to the true distribution. In particular, if no block of that length has occurred yet, clearly no intelligent estimate can be given. For this reason we will estimate only along stopping times $\lambda_1, \lambda_2, \dots$ and our main positive result is that there is a sequence of universally defined stopping times $\lambda_n$ with density 1 and estimators $h_n(X_0, X_1, \dots, X_{\lambda_n})$ which are almost surely converging to $\mu_{\tau(X_0, X_1, \dots, X_{\lambda_n})}$. (For further reading on estimation along stopping times see [12, 15, 16, 18, 21, 22].) We also will define estimators $\hat{p}_l(X_0, X_1, \dots, X_{\lambda_n})$ which are almost surely converging in the variation metric to the conditional distribution of the residual waiting time. These results will require a suitable higher moment condition on the $\{p_k\}$ distribution. These estimators are simply the averages of what we observe in a piece of the data segment $X_{\kappa_n}, \dots, X_{\lambda_n}$ where $\kappa_n$ is chosen so that there is a large fixed number of occurrences of the relevant pattern. The reason for these stopping times $\lambda_n$ is that we want to estimate only at those times when we feel that we have enough data.

Another kind of result may be obtained without a higher moment condition. Namely, there is a sequence of estimators $\tilde{h}_n$ and $\tilde{p}_n$ such that for any renewal process and almost every sequence of observations $X_0, X_1, \dots$ there is a sequence of density 1 of $n$'s $D$, which depend on the observed sequence of $X_i$ along which these estimators converge to the $\mu_\tau$ and conditional distributions of residual waiting times. The difference is that now we are unable to determine what these sequences are by finite observations.

On the other hand, for stopping times of density 1 we will show that no such result is possible in general, that is, without higher moment assumptions. More precisely, there is no strictly increasing sequence of stopping times $\{\lambda_n\}$ with density 1, and sequence of estimators $\{h_n(X_0, \dots, X_{\lambda_n})\}$, such that for all binary classical renewal processes

$$\limsup_{n\to\infty} |h_n(X_0, \dots, X_{\lambda_n}) - \mu_{\tau(X_0, \dots, X_{\lambda_n})}| = 0 \qquad \text{almost surely.}$$

(For results of similar vein see [4, 9, 19, 20, 27].)

In spite of this negative result, without any condition on higher moments, we can find stopping times with density close to 1 along which we converge to the estimates that are possible with full knowledge of the system. That is to say, for any $\varepsilon > 0$ there exists a sequence of stopping times $\{\lambda_n^{(\varepsilon)}\}$, estimators $\{h_n^{(\varepsilon)}(X_0, \dots, X_{\lambda_n^{(\varepsilon)}})\}$ and $\{p_l^{(\varepsilon)}(X_0, \dots, X_{\lambda_n^{(\varepsilon)}})\}$ such that the density of the stopping times is greater than $1 - \varepsilon$ and almost surely these estimators converge to $\mu_\tau$ and the conditional distributions of residual waiting times, respectively.

**2. Results.** It is easiest to formally define a renewal process in terms of an underlying Markov chain. Consider a Markov chain on the state space $\{0, 1, 2, \dots\}$ with transition probabilities $p_{i,i-1} = 1$ for all $i \geq 1$ and $p_{0,i} = p_i$ a probability distribution $\pi$ on $\{0, 1, 2, \dots\}$; cf. [8], Example 12.13. This chain is positive recurrent exactly when $\sum_{i=0}^{\infty} i p_{0,i} = \mu < \infty$ and the unique stationary probability assigns mass $\frac{1}{1+\mu}$ to the state 0; cf. [7], Chapter XIII and [28], Section I.2.c. Collapsing all states $i \ge 1$ to 1 gives rise to the classical binary renewal process. Even though our primary interest is in one-sided processes, stationarity implies that there exists a two-sided process with the same statistics and we will use the two-sided version whenever it is convenient to do so.

For conciseness, we will denote $X_i^j = (X_i, \dots, X_j)$ and also use this notation for $i = -\infty$ and $j = \infty$. Our interest is in the waiting time to renewal (the state 0) given some previous observations, in particular given $X_0^n$. Recall that if the data segment $X_0^n$ does not contain a zero, the expected time to the first occurrence of a zero may be infinite; this depends on the finiteness of the second moment of $\pi$. If a zero occurs, then the expected time depends on the location of the zero and so we introduce the notation:

$$\tau(X_{-\infty}^n) = \text{the } t \ge 0 \text{ such that } X_{n-t} = 0, \quad \text{and} \quad X_i = 1 \qquad \text{for } n-t < i \le n.$$

Note that this is well defined with probability 1. If a zero occurs in $X_0^n$, then $\tau(X_{-\infty}^n)$ depends only on $X_0^n$ and so we will also write for $\tau(X_{-\infty}^n)$, $\tau(X_0^n)$ with the understanding that this is defined only if a zero occurs in $X_0^n$.

Now for the classical binary renewal process $\{X_n\}$ define $\theta_n$ as

$$\theta_n = E(\max\{0 \le k : X_i = 1 \text{ for all } n < i \le n+k\} | X_0^n).$$

(Note that $\theta_n = \frac{\sum_{k=0}^{\infty} k p_{k+\tau(X_0,\dots,X_n)}}{\sum_{k=\tau(X_0,\dots,X_n)}^{\infty} p_k}$ as soon as there is at least one zero in $X_0^n$. As we have already mentioned, if no zero occurs, then it might happen that $\theta_n = \infty$.) For a family of processes $\{X_n^{(j)}\}$ we use the notation $\theta_n^{(j)}$. Our goal is to estimate both $\theta_n$ and the distribution of the time to renewal given $X_0^n$ but without prior knowledge of the distribution function of the process.

Define $\psi$ as the position of the first zero, that is,

$$\psi = \min\{t \ge 0 : X_t = 0\}.$$

Let $0 < \gamma < 1$ be arbitrary. First define the stopping times $\lambda_n$ as $\lambda_0 = \psi$ and for $n \ge 1$,

$$\lambda_n = \min\{k > \lambda_{n-1} : |\{\psi \le i < k : \tau(X_0^i) = \tau(X_0^k)\}| \ge k^{1-\gamma}\}.$$

These are the successive times $i$ when the value $t = \tau(X_0^i)$ has occurred previously enough times so that we can safely estimate the residual renewal time by empirical distributions derived from observations already made. We also need to fix $\kappa_n$ as the index where reading backward from $X_{\lambda_n}$ will have seen for the first time $\ge \lambda_n^{1-\gamma}$ occurrences of an $i$ with $\tau(X_0^i) = \tau(X_0^{\lambda_n})$. Formally put

$$\kappa_n = \max\{K : |\{K \le k < \lambda_n : \tau(X_0^k) = \tau(X_0^{\lambda_n})\}| = \lceil \lambda_n^{1-\gamma} \rceil\}.$$

Define $\sigma_i$ as the length of runs of 1's starting at position $i$. Formally put

$$\sigma_i = \max\{0 \le l : X_j = 1 \text{ for } i < j \le i + l\}.$$

For $n > 0$ define our estimator $h_n(X_0, \dots, X_{\lambda_n})$ at time $\lambda_n$ as

$$h_n(X_0, \dots, X_{\lambda_n}) = \frac{1}{\lceil (\lambda_n)^{1-\gamma} \rceil} \sum_{i=\kappa_n}^{\lambda_n - 1} I_{\{\tau(X_0^i) = \tau(X_0^{\lambda_n})\}} \sigma_i.$$

[Notice that the role of $\kappa_n$ is rather technical. It ensures that we take into consideration exactly $\lceil (\lambda_n)^{1-\gamma} \rceil$ pieces of occurrences.] The above formula is simply the average of the residual waiting times that we have already observed in the data segment $X_{\kappa_n}^{\lambda_n}$ when we were at the same value of $\tau$ as we see at time $\lambda_n$. In a similar fashion we can define the average of the number of times that the residual waiting time assumed a fixed value. Namely, define $\hat{p}_l(X_0, \dots, X_{\lambda_n})$ for each $l$ as

$$\hat{p}_l(X_0, \dots, X_{\lambda_n}) = \frac{1}{\lceil (\lambda_n)^{1-\gamma} \rceil} \sum_{i=\kappa_n}^{\lambda_n - 1} I_{\{\tau(X_0^i) = \tau(X_0^{\lambda_n}), \sigma_i = l\}}.$$

Note that $\hat{p}_l(X_0, \dots, X_{\lambda_n})$ is a probability distribution on the nonnegative integers.

THEOREM 1. *Assume* $\sum_{k=0}^{\infty} k^{\alpha+1} p_k < \infty$ *for some* $\alpha > 2$. *Let* $0 < \gamma < \min(1 - 2/\alpha, 1/3)$. *Then for the stopping times* $\lambda_n$ *and the estimator* $h_n(X_0, \dots, X_{\lambda_n})$, $\hat{p}_l(X_0, \dots, X_{\lambda_n})$ *defined above, almost surely,*

$$\lim_{n \to \infty} \frac{\lambda_n}{n} = 1, \tag{1}$$

$$\lim_{n \to \infty} |h_n(X_0, \dots, X_{\lambda_n}) - \theta_{\lambda_n}| = 0 \tag{2}$$

*and*

$$\lim_{n \to \infty} \sum_{l=0}^{\infty} \left| \hat{p}_l(X_0, \dots, X_{\lambda_n}) - \frac{p_{l + \tau(X_0^{\lambda_n})}}{\sum_{i=\tau(X_0^{\lambda_n})}^{\infty} p_i} \right| = 0. \tag{3}$$

Note that $\hat{p}_l(X_0, \dots, X_{\lambda_n})$, $h_n$ and $\lambda_n$ depend on $\gamma$ and so on $\alpha$.

In order to reduce our assumption from $\alpha > 2$ to $\alpha > 1$ a slightly more involved scheme of stopping times is needed.

Let $0 < \gamma < 1$ be arbitrary. First define the stopping times $\lambda_n^*$ as $\lambda_0^* = \psi$ and for $n \ge 1$,

$$\lambda_n^* = \min\{t > \lambda_{n-1}^* : \exists \psi < i < \log t \text{ such that } \tau(X_0^i) = \tau(X_0^t)$$
$$\text{and } |\{\log t \le j < 2^{\lfloor \log t \rfloor} : \tau(X_0^j) = \tau(X_0^t)\}| \ge 2^{\lfloor \log t \rfloor (1-\gamma)}\}.$$

(Note that all logarithms are to the base 2.) Put

$$\kappa_n^* = \min\{K : |\{\lfloor \log \lambda_n^* \rfloor < j \le K : \tau(X_0^j) = \tau(X_0^{\lambda_n^*})\}| = \lceil 2^{\lfloor \log \lambda_n^* \rfloor (1-\gamma)} \rceil\}.$$

Note that $\kappa_n^* < 2^{\lfloor \log \lambda_n^* \rfloor}$. For $n > 0$ define our estimator $h_n^*(X_0, \dots, X_{\lambda_n^*})$ at time $\lambda_n^*$ as

$$h_n^*(X_0, \dots, X_{\lambda_n^*}) = \frac{1}{\lceil 2^{\lfloor \log \lambda_n^* \rfloor (1-\gamma)} \rceil} \sum_{i=\lfloor \log \lambda_n^* \rfloor + 1}^{\kappa_n^*} I_{\{\tau(X_0^i) = \tau(X_0^{\lambda_n^*})\}} \sigma_i.$$

(Notice that $\kappa_n^*$ ensures that we take into consideration exactly $\lceil 2^{\lfloor \log \lambda_n^* \rfloor (1-\gamma)} \rceil$ pieces of occurrences.) The above formula is simply the average of the residual waiting times that we have already observed in the data segment $X_{\lfloor \log \lambda_n^* \rfloor + 1}^{\kappa_n^*}$ when we were at the same value of $\tau$ as we see at time $\lambda_n^*$. Note that $h_n^*(X_0, \dots, X_{\lambda_n^*})$ is by far not as efficient as $h_n(X_0, \dots, X_{\lambda_n})$ since as long as $2^m \le \lambda_n^* < 2^{m+1}$ the estimator $h_n^*(X_0, \dots, X_{\lambda_n^*})$ is not refreshed. Keeping the same estimate for many values of $n$ enables us to use weaker moment assumptions since the number of unfavorable events that we have to consider is reduced.

In a similar fashion we can define the average of the number of times that the residual waiting time assumed a fixed value. Namely, define $\hat{p}_l^*(X_0, \dots, X_{\lambda_n^*})$ for each $l$ as

$$\hat{p}_l^*(X_0, \dots, X_{\lambda_n^*}) = \frac{1}{\lceil 2^{\lfloor \log \lambda_n^* \rfloor (1-\gamma)} \rceil} \sum_{i=\lfloor \log \lambda_n^* \rfloor + 1}^{\kappa_n^*} I_{\{\tau(X_0^i) = \tau(X_0^{\lambda_n^*}), \sigma_i = l\}}.$$

Note that $\hat{p}_l^*(X_0, \dots, X_{\lambda_n^*})$ is a probability distribution on the nonnegative integers.

THEOREM 2. *Assume* $\sum_{k=0}^{\infty} k^{\alpha+1} p_k < \infty$ *for some* $\alpha > 1$. *Let* $0 < \gamma < 1/3$. *Then for the stopping times* $\lambda_n^*$ *and the estimator* $h_n^*(X_0, \dots, X_{\lambda_n^*})$, $\hat{p}_l^*(X_0, \dots, X_{\lambda_n^*})$ *defined above, almost surely,*

$$\lim_{n\to\infty} \frac{\lambda_n^*}{n} = 1, \tag{4}$$

$$\lim_{n\to\infty} |h_n^*(X_0, \dots, X_{\lambda_n^*}) - \theta_{\lambda_n^*}| = 0 \tag{5}$$

*and*

$$\lim_{n\to\infty} \sum_{l=0}^{\infty} \left| \hat{p}_l^*(X_0, \dots, X_{\lambda_n^*}) - \frac{p_{l+\tau(X_0^{\lambda_n^*})}}{\sum_{i=\tau(X_0^{\lambda_n^*})}^{\infty} p_i} \right| = 0. \tag{6}$$

Note that neither $h_n^*$, $\hat{p}_l^*(X_0, \dots, X_{\lambda_n^*})$ nor $\lambda_n^*$ depend on $\alpha$.

The main point in the above theorems is that we eventually know when the error is small. If we do not want to know this, then the moment condition can be dropped as is exhibited in the next theorem.

Define the estimator $\tilde{h}_n(X_0^n)$ as

$$\tilde{h}_n(X_0^n) = \frac{\sum_{i=\psi}^{n-1} \sigma_i I_{\{\tau(X_0^i)=\tau(X_0^n)\}}}{|\{\psi \le i \le n-1 : \tau(X_0^i) = \tau(X_0^n)\}|}.$$

This is just the average of values of $\theta_i$ for the data segment $X_0^n$ for those indices $i$ for which $\tau(X_0^i) = \tau(X_0^n)$.

Define also $\tilde{p}_l(X_0^n)$ as

$$\tilde{p}_l(X_0^n) = \frac{\sum_{i=\psi}^{n-1} I_{\{\tau(X_0^i)=\tau(X_0^n),\sigma_i=l\}}}{|\{\psi \le i \le n-1 : \tau(X_0^i) = \tau(X_0^n)\}|}.$$

THEOREM 3. *For any binary renewal process $\{X_n\}$, and almost every sequence of observations $X_0^\infty$, there is a set of indices $D(X_0^\infty) \subset \{0, 1, \ldots\}$ such that $\lim_{n\to\infty} \frac{|D(X_0^\infty)\cap\{0,1,\ldots,n\}|}{n+1} = 1$ and*

$$\lim_{n\in D(X_0^\infty), n\to\infty} |\tilde{h}_n(X_0, \ldots, X_n) - \theta_n| = 0 \tag{7}$$

*and*

$$\lim_{n\in D(X_0^\infty), n\to\infty} \sum_{l=0}^{\infty} \left| \tilde{p}_l(X_0^n) - \frac{p_{l+\tau(X_0^n)}}{\sum_{i=\tau(X_0^n)}^{\infty} p_i} \right| = 0. \tag{8}$$

However, for stopping times we need some restrictions to achieve consistency on density 1 as is showed in the next theorem.

THEOREM 4. *For any strictly increasing sequence of stopping times $\{\lambda_n\}$ and sequence of estimators $\{h_n(X_0, \ldots, X_{\lambda_n})\}$, such that for all binary classical renewal processes $\lim_{n\to\infty} \frac{\lambda_n}{n} = 1$ almost surely, there exists a binary classical renewal process such that*

$$P\left( \limsup_{n\to\infty} |h_n(X_0, \ldots, X_{\lambda_n}) - \theta_{\lambda_n}| > 0 \right) > 0.$$

We do not know if a similar result can be formulated for the estimation of the distribution of the residual waiting times in total variation.

Finally, if one merely intends to predict along a stopping time with density greater than $1-\varepsilon$ for some fixed $\varepsilon > 0$, then no condition on higher moments at all is required as it is stated in the next theorem. Let $\lambda_0^{(\varepsilon)} = \psi$ and for $n \ge 1$ define

$$\lambda_n^{(\varepsilon)} = \min\{t > \lambda_{n-1}^{(\varepsilon)} : |\{\psi \le i \le t : \tau(X_0^i) < \tau(X_0^t)\}| \le t(1-\varepsilon/2)\}.$$

This sequence of stopping times is designed so that eventually we only stop when $\tau(X_0^i)$ takes values bounded by some finite $L$. The point is that if $L$ is large enough, then eventually the density of times $i$ when $\tau(X_0^i) < L$ will be greater than $(1-\varepsilon/2)$ so that our stopping times will choose only moments that are less than $L$. On the other hand, the fact that eventually the $\tau(X_0^i) < L$'s will enable us to prove the convergence of the empirical estimators by a direct application of the ergodic theorem.

Define the estimator $h_n^{(\varepsilon)}(X_0,\dots,X_{\lambda_n^{(\varepsilon)}})$ at time $\lambda_n^{(\varepsilon)}$ as

$$h_n^{(\varepsilon)}(X_0,\dots,X_{\lambda_n^{(\varepsilon)}}) = \frac{\sum_{i=\psi}^{\lambda_n^{(\varepsilon)}-1} I_{\{\tau(X_0^i)=\tau(X_0^{\lambda_n^{(\varepsilon)}})\}}\sigma_i}{|\{\psi \le i < \lambda_n^{(\varepsilon)} : \tau(X_0^i) = \tau(X_0^{\lambda_n^{(\varepsilon)}})\}|}.$$

Also define

$$p_l^{(\varepsilon)}(X_0,\dots,X_{\lambda_n^{(\varepsilon)}}) = \frac{\sum_{i=\psi}^{\lambda_n^{(\varepsilon)}-1} I_{\{\tau(X_0^i)=\tau(X_0^{\lambda_n^{(\varepsilon)}}),\sigma_i=l\}}}{|\{\psi \le i < \lambda_n^{(\varepsilon)} : \tau(X_0^i) = \tau(X_0^{\lambda_n^{(\varepsilon)}})\}|}.$$

THEOREM 5. *For the stopping times $\lambda_n^{(\varepsilon)}$ and estimator $h_n^{(\varepsilon)}(X_0,\dots,X_{\lambda_n^{(\varepsilon)}})$ defined above, almost surely,*

$$\liminf_{n\to\infty} \frac{n}{\lambda_n^{(\varepsilon)}} > 1-\varepsilon,$$

$$\limsup_{n\to\infty} |h_n^{(\varepsilon)}(X_0,\dots,X_{\lambda_n^{(\varepsilon)}}) - \theta_{\lambda_n^{(\varepsilon)}}| = 0$$

*and*

$$\limsup_{n\to\infty} \sum_{l=0}^{\infty} \left| p_l^{(\varepsilon)}(X_0,\dots,X_{\lambda_n^{(\varepsilon)}}) - \frac{p_{l+\tau(X_0^n)}}{\sum_{i=\tau(X_0^n)}^{\infty} p_i} \right| = 0.$$

**3. Proof of Theorem 1.** It is easy to see that $\lim_{n\to\infty} \frac{\lambda_n}{n} = 1$ since if a block of 1's has positive probability it will appear with that frequency which is eventually greater than $\frac{\lambda_n^{1-\gamma}}{\lambda_n}$ (which tends to zero). Formally,

$$\liminf_{n\to\infty} \frac{n}{\lambda_n} \ge \liminf_{N\to\infty} \frac{\max\{i>0 : \lambda_i \le N\}}{N}.$$

Thus we have to see why the density of times when we stop and estimate tends to 1. Since the cutoff $\frac{\lambda_n^{1-\gamma}}{\lambda_n}$ tends to zero, any positive probability event will eventually be greater than it and so for any bounded $K$ we will have

$$\liminf_{N\to\infty} \frac{\max\{i>0 : \tau(X_0^{\lambda_i}) < K, \lambda_i \le N\}}{N} = P(\tau(X_{-\infty}^0) < K).$$

As $K$ tends to $\infty$ this last expression tends to 1 and thus

$$1 \geq \limsup_{n\to\infty} \frac{n}{\lambda_n} \geq \liminf_{n\to\infty} \frac{n}{\lambda_n} \geq P(\tau(X^0_{-\infty}) < \infty) = 1.$$

This establishes (1).

The usual proof of the weak law of large numbers for independent and identically random variables $\{Z_n\}$ with a second moment uses Chebyshev's inequality $P(|\sum_{i=1}^n (Z_i - EZ_i)| \geq n\varepsilon) \leq \frac{1}{n\varepsilon^2} E((Z_1 - EZ_1)^2)$. We will need a sharpening of this for random variables with an $\alpha$th moment for $\alpha > 2$.

It will be convenient to extend our process, as we may to the past, and establish first an inequality for an estimator based on an unlimited past. For a given fixed $k$, for $i \geq 0$ define $j_i^k$ as the $i$th occurrence of $\tau(X^k_{-\infty})$ (reading backward) from position $k$, that is,

$$j_i^k = \max\{j \leq k : |\{j \leq l < k : \tau(X^l_{-\infty}) = \tau(X^k_{-\infty})\}| = i\}.$$

Now for $i \geq 0$ define

$$Z_i^{(k)} = \sigma_{j_i^k}.$$

Clearly $Z_i^{(k)}$ are conditionally independent and identically distributed given $\tau(X^k_{-\infty}) = L$. Apply Markov inequality and Theorem 2.10 of Petrov [26] to get that

$$P\left(\left|\frac{\sum_{i=1}^{\lceil k^{1-\gamma}\rceil} Z_i^{(k)}}{\lceil k^{1-\gamma}\rceil} - \frac{\sum_{h=0}^\infty h p_{h+L}}{\sum_{h=L}^\infty p_h}\right| > \varepsilon \;\middle|\; \tau(X^k_{-\infty}) = L\right)$$
$$\leq \frac{2C(\alpha)}{\varepsilon^\alpha k^{(1-\gamma)\alpha/2}} \frac{\sum_{h=0}^\infty h^\alpha p_{h+L}}{\sum_{h=L}^\infty p_h}$$

where $C(\alpha)$ depends only on $\alpha$. [Notice that $\frac{E(|Z_0^{(k)}|^\alpha | \tau(X^k_{-\infty}) = L) = \sum_{h=0}^\infty h^\alpha p_{h+L}}{\sum_{h=L}^\infty p_h}$.]

Multiply both sides of the last inequality by $P(\tau(X^k_{-\infty}) = L) = \frac{1}{1+\sum_{h=0}^\infty h p_h} \sum_{h=L}^\infty p_h$ (note that by Kac's theorem $P(X_{k-L} = 0) = \frac{1}{1+\sum_{h=0}^\infty h p_h}$; cf. [7], Chapter XIII and [28], Section I.2.c) and sum over $L$. It is easy to see that

$$\sum_{L=0}^\infty \frac{\sum_{h=0}^\infty h^\alpha p_{h+L}}{\sum_{h=L}^\infty p_h} \frac{\sum_{h=L}^\infty p_h}{1 + \sum_{h=0}^\infty h p_h} \leq \frac{\sum_{h=0}^\infty h^{\alpha+1} p_h}{1 + \sum_{h=0}^\infty h p_h}$$

and we get the following estimate:

$$P\left(\left|\frac{\sum_{i=1}^{\lceil k^{1-\gamma}\rceil} Z_i^{(k)}}{\lceil k^{1-\gamma}\rceil} - \frac{\sum_{h=0}^\infty h p_{h+\tau(X^k_{-\infty})}}{\sum_{h=\tau(X^k_{-\infty})}^\infty p_h}\right| > \varepsilon\right) \leq \frac{2C(\alpha)}{\varepsilon^\alpha k^{(1-\gamma)\alpha/2}} \frac{\sum_{h=0}^\infty h^{\alpha+1} p_h}{1 + \sum_{h=0}^\infty h p_h}.$$

Applying the Borel–Cantelli lemma [by assumption $\frac{(1-\gamma)\alpha}{2} > 1$] one gets that

$$\left|\frac{\sum_{i=1}^{\lceil k^{1-\gamma}\rceil} Z_i^{(k)}}{\lceil k^{1-\gamma}\rceil} - \frac{\sum_{h=0}^{\infty} h p_{h+\tau(X_{-\infty}^k)}}{\sum_{h=\tau(X_{-\infty}^k)}^{\infty} p_h}\right| < \varepsilon$$

eventually almost surely. Particularly, on the subsequence $\lambda_n$,

$$\left|\frac{\sum_{i=1}^{\lceil(\lambda_n)^{1-\gamma}\rceil} Z_i^{(\lambda_n)}}{\lceil(\lambda_n)^{1-\gamma}\rceil} - \frac{\sum_{h=0}^{\infty} h p_{h+\tau(X_{-\infty}^{\lambda_n})}}{\sum_{h=\tau(X_{-\infty}^{\lambda_n})}^{\infty} p_h}\right| < \varepsilon$$

eventually almost surely. Since $\tau(X_{-\infty}^{\lambda_n}) = \tau(X_0^{\lambda_n})$, $h_n(X_0^{\lambda_n}) = \frac{\sum_{i=1}^{\lceil(\lambda_n)^{1-\gamma}\rceil} Z_i^{(\lambda_n)}}{\lceil(\lambda_n)^{1-\gamma}\rceil}$ and $\theta_{\lambda_n} = \frac{\sum_{h=0}^{\infty} h p_{h+\tau(X_0^{\lambda_n})}}{\sum_{h=\tau(X_0^{\lambda_n})}^{\infty} p_h}$, we get that

$$|h_n(X_0,\dots,X_{\lambda_n}) - \theta_{\lambda_n}| < \varepsilon$$

eventually almost surely which since $\varepsilon$ was arbitrary gives (2).

For (3) observe that

$$\begin{aligned}
&\sum_{l=0}^{\infty}\left|\hat{p}_l(X_0^{\lambda_n}) - \frac{p_{l+\tau(X_0^{\lambda_n})}}{\sum_{i=\tau(X_0^{\lambda_n})}^{\infty} p_i}\right| \\
&\quad = \sum_{l=0}^{\lceil(\lambda_n)^{\gamma}\log\lambda_n\rceil-1}\left|\hat{p}_l(X_0^{\lambda_n}) - \frac{p_{l+\tau(X_0^{\lambda_n})}}{\sum_{i=\tau(X_0^{\lambda_n})}^{\infty} p_i}\right| \\
&\qquad + \sum_{l=\lceil(\lambda_n)^{\gamma}\log\lambda_n\rceil}^{\infty}\left|\hat{p}_l(X_0^{\lambda_n}) - \frac{p_{l+\tau(X_0^{\lambda_n})}}{\sum_{i=\tau(X_0^{\lambda_n})}^{\infty} p_i}\right| \\
&\quad = A_n + B_n.
\end{aligned}$$

First we deal with the first term. We will use finite sums of exponential bounds in order to bound it. Now define

$$Z_{i,l}^{(k)} = I_{\{\sigma_{j_i^k}=l\}}.$$

Clearly $Z_{i,l}^{(k)}$ are conditionally independent and identically distributed given $\tau(X_{-\infty}^k) = L$. Apply Hoeffding's inequality to get that

$$\begin{aligned}
&P\left(\left|\frac{\sum_{i=1}^{\lceil k^{1-\gamma}\rceil} Z_{i,l}^{(k)}}{\lceil k^{1-\gamma}\rceil} - \frac{p_{l+L}}{\sum_{h=L}^{\infty} p_h}\right| > k^{-\gamma}(\log k)^{-2}\middle|\tau(X_{-\infty}^k) = L\right) \\
&\quad \le e^{-k^{1-\gamma}/(2k^{2\gamma}(\log k)^4)}.
\end{aligned}$$

After integrating both sides with respect to the conditioning, and using the sum bound on the events for $0 \le l \le \lceil k^\gamma \log k\rceil - 1$, we get

$$P\left(\max_{0\le l\le \lceil k^\gamma \log k\rceil-1}\left|\frac{\sum_{i=1}^{\lceil k^{1-\gamma}\rceil} Z_{i,l}^{(k)}}{\lceil k^{1-\gamma}\rceil} - \frac{p_{l+\tau(X_{-\infty}^k)}}{\sum_{h=\tau(X_{-\infty}^k)}^\infty p_h}\right| > k^{-\gamma}(\log k)^{-2}\right)$$
$$\le \lceil k^\gamma \log k\rceil e^{-k^{1-\gamma}/(2k^{2\gamma}(\log k)^4)}$$

which is summable (by assumption $\gamma < \frac{1}{3}$ ) and so by the Borel–Cantelli lemma,

$$\max_{0\le l\le \lceil k^\gamma \log k\rceil-1}\left|\frac{\sum_{i=1}^{\lceil k^{1-\gamma}\rceil} Z_{i,l}^{(k)}}{\lceil k^{1-\gamma}\rceil} - \frac{p_{l+\tau(X_{-\infty}^k)}}{\sum_{h=\tau(X_{-\infty}^k)}^\infty p_h}\right| \le k^{-\gamma}(\log k)^{-2}$$

eventually almost surely. Particularly, on the subsequence $\lambda_n$,

$$\sum_{l=0}^{\lceil(\lambda_n)^\gamma \log \lambda_n\rceil-1}\left|\frac{\sum_{i=1}^{\lceil(\lambda_n)^{1-\gamma}\rceil} Z_{i,l}^{(\lambda_n)}}{\lceil(\lambda_n)^{1-\gamma}\rceil} - \frac{p_{l+\tau(X_{-\infty}^{\lambda_n})}}{\sum_{h=\tau(X_{-\infty}^{\lambda_n})}^\infty p_h}\right|$$
$$\le \lceil(\lambda_n)^\gamma(\log\lambda_n)\rceil {\lambda_n}^{-\gamma}(\log\lambda_n)^{-2}$$

eventually almost surely. Observe that $\tau(X_{-\infty}^{\lambda_n}) = \tau(X_0^{\lambda_n})$ and $\hat{p}_l(X_0^{\lambda_n}) = \frac{\sum_{i=1}^{\lceil(\lambda_n)^{1-\gamma}\rceil} Z_{i,l}^{(\lambda_n)}}{\lceil(\lambda_n)^{1-\gamma}\rceil}$ and so we get that

$$A_n = \sum_{l=0}^{\lceil(\lambda_n)^\gamma \log \lambda_n\rceil-1}\left|\hat{p}_l(X_0^{\lambda_n}) - \frac{p_{l+\tau(X_0^{\lambda_n})}}{\sum_{i=\tau(X_0^{\lambda_n})}^\infty p_i}\right| \le \frac{2}{\log\lambda_n} \tag{9}$$

eventually almost surely. We have to prove that $B_n \to 0$ almost surely. Note that by the Markov inequality, given $\tau(X_0^k) = L$, for $L < k$,

$$\frac{\sum_{l=\lceil\mu_L \log k\rceil}^\infty p_{l+L}}{\sum_{l=L}^\infty p_l} \le \frac{1}{\log k} \tag{10}$$

where $\mu_L = \sum_{i=L}^\infty (i-L)p_i / \sum_{i=L}^\infty p_i$.

Now observe that almost surely for sufficiently large $n$,

$$\mu_{\tau(X_0^{\lambda_n})} \le (\lambda_n)^\gamma. \tag{11}$$

Indeed

$$h_n(X_0^{\lambda_n}) = \frac{1}{\lceil(\lambda_n)^{1-\gamma}\rceil}\sum_{i=\kappa_n}^{\lambda_n-1} I_{\{\tau(X_0^i)=\tau(X_0^{\lambda_n})\}}\sigma_{i,\lambda_n} \le \frac{\lambda_n - \lceil(\lambda_n)^{1-\gamma}\rceil}{\lceil(\lambda_n)^{1-\gamma}\rceil} \le (\lambda_n)^\gamma - 1$$

[in the data segment $X_0^{\lambda_n}$ there are at least $\lceil(\lambda_n)^{1-\gamma}\rceil$ zeros] and we have already proved that $h_n(X_0^{\lambda_n}) - \mu_{\tau(X_0^{\lambda_n})} \to 0$.

Now, apply (11), (10) and the upper bound on $A_n$ in (9) in order to get

$$\begin{aligned}
B_n &\le \sum_{l=\lceil(\lambda_n)^\gamma \log\lambda_n\rceil}^{\infty} \hat{p}_l(X_0^{\lambda_n}) + \sum_{l=\lceil(\lambda_n)^\gamma \log\lambda_n\rceil}^{\infty} \frac{p_{l+\tau(X_0^{\lambda_n})}}{\sum_{i=\tau(X_0^{\lambda_n})}^{\infty} p_i} \\
&\le 1 - \sum_{l=0}^{\lceil(\lambda_n)^\gamma \log\lambda_n\rceil-1} \hat{p}_l(X_0^{\lambda_n}) + \frac{2}{\log\lambda_n} \\
&\le 1 - \sum_{l=0}^{\lceil(\lambda_n)^\gamma \log\lambda_n\rceil-1} \frac{p_{l+\tau(X_0^{\lambda_n})}}{\sum_{i=\tau(X_0^{\lambda_n})}^{\infty} p_i} \\
&\quad + \sum_{l=0}^{\lceil(\lambda_n)^\gamma \log\lambda_n\rceil-1} \left| \hat{p}_l(X_0^{\lambda_n}) - \frac{p_{l+\tau(X_0^{\lambda_n})}}{\sum_{i=\tau(X_0^{\lambda_n})}^{\infty} p_i} \right| + \frac{2}{\log\lambda_n} \\
&\le \frac{6}{\log n}
\end{aligned}$$

eventually almost surely, and so $B_n \to 0$ almost surely. The proof of Theorem 1 is complete.

**4. Proof of Theorem 2.** The proof is similar to that of Theorem 1 but with a number of changes required to deal with the weaker hypothesis. It is easy to see that $\lim_{n\to\infty} \frac{\lambda_n^*}{n} = 1$ since if a block of 1's has positive probability it will appear with that frequency which is eventually greater than $\frac{2^{\lfloor \log \lambda_n^* \rfloor (1-\gamma)}}{\lambda_n^*}$ (which tends to zero). Formally,

$$\begin{aligned}
\liminf_{n\to\infty} \frac{n}{\lambda_n^*} &\ge \liminf_{N\to\infty} \frac{\max\{i>0 : \lambda_i^* \le N\}}{N} \\
&\ge \liminf_{N\to\infty} \frac{\max\{i>0 : \tau(X_0^{\lambda_i^*}) < K, \lambda_i^* \le N\}}{N} \\
&= P(\tau(X_{-\infty}^0) < K)
\end{aligned}$$

for arbitrary large $K$. Thus

$$1 \ge \limsup_{n\to\infty} \frac{n}{\lambda_n^*} \ge \liminf_{n\to\infty} \frac{n}{\lambda_n^*} \ge P(\tau(X_{-\infty}^0) < \infty) = 1.$$

Let $k < m$ be fixed. Define $j_0^{(k,m)} = m$ and for $i \ge 0$ let $j_{i+1}^{(k,m)}$ denote the $(i+1)$st occurrence of $\tau(X_{-\infty}^k)$ (reading forward, starting at position $m$), that is,

$$j_{i+1}^{(k,m)} = \min\{t > j_i^{(k,m)} : \tau(X_{-\infty}^t) = \tau(X_{-\infty}^k)\}.$$

Now for $i \ge 1$ define

$$Z_i^{(k,m)} = \sigma_{j_i^{(k,m)}}.$$

Clearly $Z_i^{(k,m)}$ are conditionally independent and identically distributed given $\tau(X_{-\infty}^k) = L$. For $1 < \alpha \le 2$ apply Markov inequality and Theorem 2 of von Bahr and Essen in [5] to get that

$$P\left( \left| \frac{\sum_{i=1}^{\lceil (2^m)^{1-\gamma} \rceil} Z_i^{(k,m)}}{\lceil (2^m)^{1-\gamma} \rceil} - \frac{\sum_{h=0}^{\infty} h p_{h+L}}{\sum_{h=L}^{\infty} p_h} \right| > \varepsilon \middle| \tau(X_{-\infty}^k) = L \right)$$
$$\le \frac{10}{\varepsilon^\alpha (2^m)^{(1-\gamma)(\alpha-1)}} \frac{\sum_{h=0}^{\infty} h^\alpha p_{h+L}}{\sum_{h=L}^{\infty} p_h}.$$

[Notice that $E(|Z_1^{(k,m)}|^\alpha | \tau(X_{-\infty}^k) = L) = \frac{\sum_{h=0}^{\infty} h^\alpha p_{h+L}}{\sum_{h=L}^{\infty} p_h}$.] Multiply both sides of the last inequality by $P(\tau(X_{-\infty}^k) = L) = \frac{1}{1+\sum_{h=0}^{\infty} h p_h} \sum_{h=L}^{\infty} p_h$ (note that by Kac's theorem $P(X_{k-L} = 0) = \frac{1}{1+\sum_{h=0}^{\infty} h p_h}$; cf. [7], Chapter XIII and [28], Section I.2.c) and sum over $L$. It is easy to see that

$$\sum_{L=0}^{\infty} \frac{\sum_{h=0}^{\infty} h^\alpha p_{h+L}}{\sum_{h=L}^{\infty} p_h} \frac{\sum_{h=L}^{\infty} p_h}{1+\sum_{h=0}^{\infty} h p_h} \le \frac{\sum_{h=0}^{\infty} h^{\alpha+1} p_h}{1+\sum_{h=0}^{\infty} h p_h}$$

and we get the following estimate:

$$P\left( \left| \frac{\sum_{i=1}^{\lceil (2^m)^{1-\gamma} \rceil} Z_i^{(k,m)}}{\lceil (2^m)^{1-\gamma} \rceil} - \frac{\sum_{h=0}^{\infty} h p_{h+\tau(X_{-\infty}^k)}}{\sum_{h=\tau(X_{-\infty}^k)}^{\infty} p_h} \right| > \varepsilon \right)$$
$$\le \frac{10}{\varepsilon^\alpha (2^m)^{(1-\gamma)(\alpha-1)}} \frac{\sum_{h=0}^{\infty} h^{\alpha+1} p_h}{1+\sum_{h=0}^{\infty} h p_h}$$

and in turn

$$P\left( \max_{0 \le k \le m-1} \left| \frac{\sum_{i=1}^{\lceil (2^m)^{1-\gamma} \rceil} Z_i^{(k,m)}}{\lceil (2^m)^{1-\gamma} \rceil} - \frac{\sum_{h=0}^{\infty} h p_{h+\tau(X_{-\infty}^k)}}{\sum_{h=\tau(X_{-\infty}^k)}^{\infty} p_h} \right| > \varepsilon \right)$$
$$\le \frac{10m}{\varepsilon^\alpha (2^m)^{(1-\gamma)(\alpha-1)}} \frac{\sum_{h=0}^{\infty} h^{\alpha+1} p_h}{1+\sum_{h=0}^{\infty} h p_h}$$

and the right-hand side is summable. For $\alpha > 2$ apply Markov inequality and Theorem 2.10 of Petrov [26] to get that

$$P\left( \left| \frac{\sum_{i=1}^{\lceil (2^m)^{1-\gamma} \rceil} Z_i^{(k,m)}}{\lceil (2^m)^{1-\gamma} \rceil} - \frac{\sum_{h=0}^{\infty} h p_{h+L}}{\sum_{h=L}^{\infty} p_h} \right| > \varepsilon \middle| \tau(X_{-\infty}^k) = L \right)$$
$$\le \frac{2C(\alpha)}{\varepsilon^\alpha 2^{m(1-\gamma)\alpha/2}} \frac{\sum_{h=0}^{\infty} h^\alpha p_{h+L}}{\sum_{h=L}^{\infty} p_h}$$

where $C(\alpha)$ depends only on $\alpha$. Integrating both sides, just as in the previous case above, we get

$$P\left(\left|\frac{\sum_{i=1}^{\lceil (2^m)^{1-\gamma}\rceil} Z_i^{(k,m)}}{\lceil (2^m)^{1-\gamma}\rceil} - \frac{\sum_{h=0}^{\infty} hp_{h+\tau(X_{-\infty}^k)}}{\sum_{h=\tau(X_{-\infty}^k)}^{\infty} p_h}\right| > \varepsilon\right)$$
$$\leq \frac{2C(\alpha)}{\varepsilon^{\alpha} 2^{m(1-\gamma)\alpha/2}} \frac{\sum_{h=0}^{\infty} h^{\alpha+1} p_h}{1+\sum_{h=0}^{\infty} hp_h}$$

and in turn

$$P\left(\max_{0\leq k\leq m-1}\left|\frac{\sum_{i=1}^{\lceil (2^m)^{1-\gamma}\rceil} Z_i^{(k,m)}}{\lceil (2^m)^{1-\gamma}\rceil} - \frac{\sum_{h=0}^{\infty} hp_{h+\tau(X_{-\infty}^k)}}{\sum_{h=\tau(X_{-\infty}^k)}^{\infty} p_h}\right| > \varepsilon\right)$$
$$\leq \frac{2mC(\alpha)}{\varepsilon^{\alpha} 2^{m(1-\gamma)\alpha/2}} \frac{\sum_{h=0}^{\infty} h^{\alpha+1} p_h}{1+\sum_{h=0}^{\infty} hp_h}$$

and the right-hand side is summable. Applying the Borel–Cantelli lemma in both cases one gets that

$$\max_{0\leq k\leq m-1}\left|\frac{\sum_{i=1}^{\lceil (2^m)^{1-\gamma}\rceil} Z_i^{(k,m)}}{\lceil (2^m)^{1-\gamma}\rceil} - \frac{\sum_{h=0}^{\infty} hp_{h+\tau(X_{-\infty}^k)}}{\sum_{h=\tau(X_{-\infty}^k)}^{\infty} p_h}\right| < \varepsilon$$

eventually almost surely. Since $2^m \leq \lambda_n^* < 2^{m+1}$ for some $m$, we get that

$$|h_n^*(X_0,\dots,X_{\lambda_n^*}) - \theta_{\lambda_n^*}| < \varepsilon$$

eventually almost surely, which since $\varepsilon$ was arbitrary gives (5). [Indeed, observe first that for $k \geq \psi$, $\tau(X_{-\infty}^k) = \tau(X_0^k)$. Now for suitable $k < \lfloor \log \lambda_n^* \rfloor$ and $m = \lfloor \log \lambda_n^* \rfloor$: $h_n^*(X_0,\dots,X_{\lambda_n^*}) = \frac{\sum_{i=1}^{\lceil (2^m)^{1-\gamma}\rceil} Z_i^{(k,m)}}{\lceil (2^m)^{1-\gamma}\rceil}$ and $\theta_{\lambda_n^*} = \frac{\sum_{h=0}^{\infty} hp_{h+\tau(X_{-\infty}^k)}}{\sum_{h=\tau(X_{-\infty}^k)}^{\infty} p_h}$.]

Now we will deal with (6). For $k < m$ define

$$Z_{i,l}^{(k,m)} = I_{\{\sigma_{j_i^{(k,m)}} = l\}}.$$

Clearly, for fixed $k < m$ and $l$, $Z_{i,l}^{(k,m)}$, $i \geq 1$, are conditionally independent and identically distributed given $\tau(X_{-\infty}^k) = L$. Apply Hoeffding's inequality to get that

$$P\left(\left|\frac{\sum_{i=1}^{\lceil (2^m)^{1-\gamma}\rceil} Z_{i,l}^{(k,m)}}{\lceil (2^m)^{1-\gamma}\rceil} - \frac{p_{l+L}}{\sum_{h=L}^{\infty} p_h}\right| > (2^m)^{-\gamma} m^{-2} \middle| \tau(X_{-\infty}^k) = L\right)$$
$$\leq e^{-(2^m)^{1-\gamma}/2(2^m)^{2\gamma} m^4}.$$

After integrating both sides with respect to the conditioning, and using the sum bound on the events for $0 \le l < \lceil (2^m)^\gamma m \rceil$, we get

$$P\left( \max_{0\le l<\lceil (2^m)^\gamma m\rceil} \left| \frac{\sum_{i=1}^{\lceil (2^m)^{1-\gamma}\rceil} Z_{i,l}^{(k,m)}}{\lceil (2^m)^{1-\gamma}\rceil} - \frac{p_{l+\tau(X_{-\infty}^k)}}{\sum_{h=\tau(X_{-\infty}^k)}^{\infty} p_h} \right| > (2^m)^{-\gamma} m^{-2} \right)$$
$$\le \lceil (2^m)^\gamma m\rceil e^{-(2^m)^{1-\gamma}/(2(2^m)^{2\gamma} m^4)}.$$

Now

$$P\left( \sum_{l=0}^{\lceil (2^m)^\gamma m\rceil -1} \left| \frac{\sum_{i=1}^{\lceil (2^m)^{1-\gamma}\rceil} Z_{i,l}^{(k,m)}}{\lceil (2^m)^{1-\gamma}\rceil} - \frac{p_{l+\tau(X_{-\infty}^k)}}{\sum_{h=\tau(X_{-\infty}^k)}^{\infty} p_h} \right| > \frac{\lceil (2^m)^\gamma m\rceil}{(2^m)^\gamma m^2} \right)$$
$$\le \lceil (2^m)^\gamma m\rceil e^{-(2^m)^{1-\gamma}/(2(2^m)^{2\gamma} m^4)}$$

and

$$P\left( \max_{0\le k<m} \sum_{l=0}^{\lceil (2^m)^\gamma m\rceil -1} \left| \frac{\sum_{i=1}^{\lceil (2^m)^{1-\gamma}\rceil} Z_{i,l}^{(k,m)}}{\lceil (2^m)^{1-\gamma}\rceil} - \frac{p_{l+\tau(X_{-\infty}^k)}}{\sum_{h=\tau(X_{-\infty}^k)}^{\infty} p_h} \right| > \frac{\lceil (2^m)^\gamma m\rceil}{(2^m)^\gamma m^2} \right)$$
$$\le \frac{m\lceil (2^m)^\gamma m\rceil}{e^{(2^m)^{1-\gamma}/(2(2^m)^{2\gamma} m^4)}},$$

which is summable and so by the Borel–Cantelli lemma,

$$\max_{0\le k<m} \sum_{l=0}^{\lceil (2^m)^\gamma m\rceil -1} \left| \frac{\sum_{i=1}^{\lceil (2^m)^{1-\gamma}\rceil} Z_{i,l}^{(k,m)}}{\lceil (2^m)^{1-\gamma}\rceil} - \frac{p_{l+\tau(X_{-\infty}^k)}}{\sum_{h=\tau(X_{-\infty}^k)}^{\infty} p_h} \right| \le \frac{\lceil (2^m)^\gamma m\rceil}{(2^m)^\gamma m^2} \le \frac{2}{m}$$

eventually almost surely. Since $2^m \le \lambda_n^* < 2^{m+1}$ for some $m$,

$$\sum_{l=0}^{\lceil 2^{\lfloor \log \lambda_n^*\rfloor \gamma} \lfloor \log \lambda_n^*\rfloor\rceil -1} \left| \hat{p}_l(X_0^{\lambda_n^*}) - \frac{p_{l+\tau(X_0^{\lambda_n^*})}}{\sum_{i=\tau(X_0^{\lambda_n^*})}^{\infty} p_i} \right| \le \frac{2}{\lfloor \log \lambda_n^*\rfloor} \tag{12}$$

eventually almost surely. [Indeed, observe first that for $k \ge \psi$, $\tau(X_{-\infty}^k) = \tau(X_0^k)$. Now for suitable $k < \lfloor \log \lambda_n^*\rfloor$ and $m = \lfloor \log \lambda_n^*\rfloor$: $\hat{p}_l(X_0^{\lambda_n^*}) = \frac{\sum_{i=1}^{\lceil (2^m)^{1-\gamma}\rceil} Z_{i,l}^{(k,m)}}{\lceil (2^m)^{1-\gamma}\rceil}$ and $\frac{p_{l+\tau(X_0^{\lambda_n^*})}}{\sum_{i=\tau(X_0^{\lambda_n^*})}^{\infty} p_i} = \frac{p_{l+\tau(X_{-\infty}^k)}}{\sum_{h=\tau(X_{-\infty}^k)}^{\infty} p_h}$.]

Observe that

$$\sum_{l=0}^{\infty} \left| \hat{p}_l(X_0^{\lambda_n^*}) - \frac{p_{l+\tau(X_0^{\lambda_n^*})}}{\sum_{i=\tau(X_0^{\lambda_n^*})}^{\infty} p_i} \right|$$

$$= \sum_{l=0}^{\lceil 2^{\lfloor \log \lambda_n^* \rfloor \gamma} \lfloor \log \lambda_n^* \rfloor \rceil - 1} \left| \hat{p}_l(X_0^{\lambda_n^*}) - \frac{p_{l+\tau(X_0^{\lambda_n^*})}}{\sum_{i=\tau(X_0^{\lambda_n^*})}^{\infty} p_i} \right|$$
$$+ \sum_{l=\lceil 2^{\lfloor \log \lambda_n^* \rfloor \gamma} \lfloor \log \lambda_n^* \rfloor \rceil}^{\infty} \left| \hat{p}_l(X_0^{\lambda_n^*}) - \frac{p_{l+\tau(X_0^{\lambda_n^*})}}{\sum_{i=\tau(X_0^{\lambda_n^*})}^{\infty} p_i} \right|$$
$$= A_n + B_n.$$

By (12), $A_n \to 0$ almost surely. We have to prove that $B_n \to 0$ almost surely. Note that by the Markov inequality, given $\tau(X_0^k) = L$, for $L < k$,

$$\frac{\sum_{l=\lceil \mu_L \lfloor \log k \rfloor \rceil}^{\infty} p_{l+L}}{\sum_{l=L}^{\infty} p_l} \le \frac{1}{\lfloor \log k \rfloor} \tag{13}$$

where $\mu_L = \sum_{i=L}^{\infty} (i-L) p_i / \sum_{i=L}^{\infty} p_i$.

Now observe that almost surely for sufficiently large $n$,

$$\mu_{\tau(X_0^{\lambda_n^*})} \le 2^{\lfloor \log \lambda_n^* \rfloor \gamma}. \tag{14}$$

Indeed

$$h_n(X_0^{\lambda_n^*}) = \frac{\sum_{i=\lfloor \log \lambda_n^* \rfloor + 1}^{\kappa_n^*} I_{\{\tau(X_0^i) = \tau(X_0^{\lambda_n^*})\}} \sigma_i}{\lceil 2^{\lfloor \log \lambda_n^* \rfloor (1-\gamma)} \rceil} \le \frac{2^{\lfloor \log \lambda_n^* \rfloor} - \lceil 2^{\lfloor \log \lambda_n^* \rfloor (1-\gamma)} \rceil}{\lceil 2^{\lfloor \log \lambda_n^* \rfloor (1-\gamma)} \rceil}$$
$$\le 2^{\lfloor \log \lambda_n^* \rfloor \gamma} - 1$$

(in the data segment $X_0^{\lambda_n^*}$ there are at least $\lceil 2^{\lfloor \log \lambda_n^* \rfloor (1-\gamma)} \rceil$ zeros) and we have already proved that $h_n(X_0^{\lambda_n^*}) - \mu_{\tau(X_0^{\lambda_n^*})} \to 0$.

Now, apply (14), (13) and the upper bound on $A_n$ in (12) in order to get

$$B_n \le \sum_{l=\lceil 2^{\lfloor \log \lambda_n^* \rfloor \gamma} \lfloor \log \lambda_n^* \rfloor \rceil}^{\infty} \hat{p}_l(X_0^{\lambda_n^*}) + \sum_{l=\lceil 2^{\lfloor \log \lambda_n^* \rfloor \gamma} \lfloor \log \lambda_n^* \rfloor \rceil}^{\infty} \frac{p_{l+\tau(X_0^{\lambda_n^*})}}{\sum_{i=\tau(X_0^{\lambda_n^*})}^{\infty} p_i}$$
$$\le 1 - \sum_{l=0}^{\lceil 2^{\lfloor \log \lambda_n^* \rfloor \gamma} \lfloor \log \lambda_n^* \rfloor \rceil - 1} \hat{p}_l(X_0^{\lambda_n^*}) + \frac{1}{\lfloor \log \lambda_n^* \rfloor}$$
$$\le 1 - \sum_{l=0}^{\lceil 2^{\lfloor \log \lambda_n^* \rfloor \gamma} \lfloor \log \lambda_n^* \rfloor \rceil - 1} \frac{p_{l+\tau(X_0^{\lambda_n^*})}}{\sum_{i=\tau(X_0^{\lambda_n^*})}^{\infty} p_i}$$
$$+ \sum_{l=0}^{\lceil 2^{\lfloor \log \lambda_n^* \rfloor \gamma} \lfloor \log \lambda_n^* \rfloor \rceil - 1} \left| \hat{p}_l(X_0^{\lambda_n^*}) - \frac{p_{l+\tau(X_0^{\lambda_n^*})}}{\sum_{i=\tau(X_0^{\lambda_n^*})}^{\infty} p_i} \right| + \frac{1}{\lfloor \log \lambda_n^* \rfloor}$$

$$\leq \sum_{l=\lceil 2^{\lfloor \log \lambda_n^* \rfloor \gamma} \lfloor \log \lambda_n^* \rfloor \rceil}^{\infty} \frac{p_{l+\tau(X_0^{\lambda_n^*})}}{\sum_{i=\tau(X_0^{\lambda_n^*})}^{\infty} p_i} + \frac{3}{\lfloor \log \lambda_n^* \rfloor}$$

$$\leq \frac{4}{\lfloor \log n \rfloor}$$

eventually almost surely, and so $B_n \to 0$ almost surely. The proof of Theorem 2 is complete.

**5. Proof of Theorem 3.** In the proof of this theorem we do not need to use explicit estimates and can rely on the ergodic theorem alone. Notice that these renewal processes are always ergodic and therefore any finite block that occurs at all with positive probability will almost surely eventually occur in the data segment $X_0^n$ with an empirical distribution which is converging to its probability. This observation yields the following for any fixed $m$:

$$\lim_{L\to\infty} \lim_{n\to\infty} \frac{1}{n} \sum_{i=0}^{n-1} I_{\{\tau(X_0^i) \leq L, |\tilde{h}_i(X_0^i) - \sum_{k=0}^{\infty} kp_{k+\tau(X_0^i)} / \sum_{k=\tau(X_0^i)}^{\infty} p_k| < 2^{-m}\}} = 1$$

almost surely. What follows is that for each $m$,

$$\lim_{n\to\infty} \frac{1}{n} \sum_{i=0}^{n-1} I_{\{|\tilde{h}_i(X_0^i) - \sum_{k=0}^{\infty} kp_{k+\tau(X_0^i)} / \sum_{k=\tau(X_0^i)}^{\infty} p_k| < 2^{-m}\}} = 1$$

almost surely.

To obtain the set $D_1$ with density 1 we will construct an auxiliary sequence of integers $N_m$ tending to infinity as follows. For a fixed realization $X_0^\infty$, let $N_0 = 0$ and for $m \geq 1$ define

$$N_m = \min\Bigg\{ n > N_{m-1} : \forall i \geq n,$$
$$\frac{1}{i+1} \sum_{j=0}^{i} I_{\{|\tilde{h}_j(X_0^j) - \sum_{k=0}^{\infty} kp_{k+\tau(X_0^j)} / \sum_{k=\tau(X_0^j)}^{\infty} p_k| < 2^{-(m+1)}\}} > 1 - 2^{-(m+1)} \Bigg\}.$$

The existence of these $N_m$'s follows once again from the ergodic theorem and since we are requiring only a countable number of conditions we may assume that these are satisfied simultaneously on a single set with probability 1. Notice that for any $i \geq N_m$ the number of indices $j$ where the error we are making is at most $2^{-(m+1)}$ is at least $j(1 - 2^{-(m+1)})$. Using this sequence define the set of indexes $D_1(X_0^\infty)$ as

$$D_1(X_0^\infty) = \bigcup_{i=1}^{\infty} \left\{ n \leq N_i : \left| \tilde{h}_n(X_0, \dots, X_n) - \frac{\sum_{k=0}^{\infty} kp_{k+\tau(X_0^n)}}{\sum_{k=\tau(X_0^n)}^{\infty} p_k} \right| < 2^{-i} \right\}.$$

By our previous observation the density of this $D_1$ will be 1, namely:

$$\lim_{n\to\infty} \frac{|D_1(X_0^\infty) \cap \{0,1,\dots,n\}|}{n+1} = 1.$$

Furthermore,

$$\lim_{n\in D(X_0^\infty), n\to\infty} \left| \tilde{h}_n(X_0,\dots,X_n) - \frac{\sum_{k=0}^\infty k p_{k+\tau(X_0^n)}}{\sum_{k=\tau(X_0^n)}^\infty p_k} \right| = 0.$$

For $\tilde{p}_l(X_0^n)$ the proof proceeds along similar lines. A set $D_2(X_0^\infty)$ is constructed with density 1 along which (8) will hold and the set $D$ in the theorem is taken to be $D_1 \cap D_2$ which has density 1. The proof of Theorem 3 is complete.

**6. Proof of Theorem 4.** Suppose that on the contrary

$$P\left( \lim_{n\to\infty} |h_n(X_0,\dots,X_{\lambda_n}) - \theta_{\lambda_n}| = 0 \right) = 1$$

for all binary classical renewal processes.

We first define an auxiliary Markov chain $\mathcal{M}^{(0)}$. Let the state space be the nonnegative integers. For $i \ge 0$ let $p_{0,i}^{(0)} = \frac{1}{2^{i+1}}$ and $p_{i+1,i}^{(0)} = 1$. Clearly, state zero is positive recurrent and since the Markov chain is irreducible this Markov chain yields a stationary and ergodic distribution. We will modify this Markov chain $\mathcal{M}^{(0)}$ in such a way that the limiting Markov chain $\mathcal{M}^{(\infty)}$ will remain stationary and ergodic.

The binary classical renewal process is defined as $X_n^{(i)} = 0$ if $M_n^{(i)} = 0$ and $X_n^{(i)} = 1$ otherwise. Let $L_0 = 0$.

Now choose $N_1$ large enough that

$$P\Bigg( \bigcup_{n=1}^\infty \Big\{ L_0 < \lambda_n < N_1, X_{\lambda_n}^{(0)} = 0, \\ \left| h_n(X_0^{(0)},\dots,X_{\lambda_n}^{(0)}) - \sum_{i=0}^\infty i p_{0,i}^{(0)} \right| < \tfrac{1}{100} \Big\} \;\Big|\; X_0^{(0)} = 0 \Bigg) > 1 - \tfrac{1}{1000}.$$

This can be done since $P(X_0^{(0)} = 0) > 0$ and $\lim_{n\to\infty} \frac{\lambda_n}{n} = 1$. Note that if $X_{\lambda_n}^{(0)} = 0$, then $\theta_{\lambda_n}^{(0)} = \sum_{i=0}^\infty i p_{0,i}^{(0)}$.

For an arbitrary $\delta_1 < 0.25 p_{0,0}^{(0)}$ (which will be specified later) let $p_{0,0}^{(1)} = p_{0,0}^{(0)} - \delta_1$ and for some $k_1 > \frac{2}{\delta_1}$, $p_{0,k_1}^{(1)} = p_{0,k_1}^{(0)} + \delta_1$. Now the change in

$$\sum_{i=0}^\infty i p_{0,i}^{(1)} - \sum_{i=0}^\infty i p_{0,i}^{(0)} = k_1 \delta_1 > 2.$$

Now choose $\delta_1$ so small such that

$$\sum_{(x_0,\dots,x_{N_1})\in\{0,1\}^{N_1}} |P(X_0^{(0)}=x_0,\dots,X_{N_1}^{(0)}=x_{N_1}|X_0^{(0)}=0)$$
$$-P(X_0^{(1)}=x_0,\dots,X_{N_1}^{(1)}=x_{N_1}|X_0^{(1)}=0)|\le \tfrac{1}{1000}.$$

In this way, for the $\{X_n^{(1)}\}$ process, for some $L_0<\lambda_n<N_1$, the estimate $h_n(X_0^{(1)},\dots,X_{\lambda_n}^{(1)})$ will be smaller than the target $\sum_{i=0}^{\infty} ip_{0,i}^{(1)}$ by at least 1 with probability $1-\frac{2}{1000}$.

For an arbitrary $N_1<L_1$ let $\sum_{i\ge L_1} p_{0,i}^{(1)}=\beta_1$. For an arbitrary $\delta_2<0.25p_{0,0}^{(1)}$ let $p_{0,0}^{(2)}=p_{0,0}^{(1)}-\delta_2$ and for some $k_2$ which will be specified later, $p_{0,k_2}^{(2)}=p_{0,k_2}^{(1)}+\delta_2$. Now the change in

$$\frac{\sum_{i=L_1}^{\infty}(i-L_1)p_{0,i}^{(2)}}{\sum_{i=L_1}^{\infty}p_{0,i}^{(2)}}-\frac{\sum_{i=L_1}^{\infty}(i-L_1)p_{0,i}^{(1)}}{\sum_{i=L_1}^{\infty}p_{0,i}^{(1)}}=\frac{k_2\delta_2}{\beta_1+\delta_2}-\frac{\delta_2\sum_{i=L_1}^{\infty}ip_{0,i}^{(1)}}{\beta_1(\beta_1+\delta_2)}$$

and in

$$\sum_{i=0}^{\infty} ip_{0,i}^{(2)}-\sum_{i=0}^{\infty} ip_{0,i}^{(1)}=k_2\delta_2.$$

Choose $N_1<L_1$ such that $3\beta_1<100^{-2}$. Now choose $N_2$ so big such that

$$P\Bigg(\bigcup_{n=1}^{\infty}\Bigg\{L_1<\lambda_n<N_2, X_{\lambda_n-L_1}^{(1)}=0, X_{\lambda_n-L_1+1}^{(1)}=\cdots=X_{\lambda_n}^{(1)}=1,$$
$$\left|h_n(X_0^{(1)},\dots,X_{\lambda_n}^{(1)})-\frac{\sum_{i=L_1}^{\infty}(i-L_1)p_{0,i}^{(1)}}{\sum_{i=L_1}^{\infty}p_{0,i}^{(1)}}\right|<\frac{1}{100}\Bigg\}\Bigg|X_0^{(1)}=0\Bigg)$$
$$>1-\left(\frac{1}{1000}\right)^2.$$

Note that if $X_{\lambda_n-L_1}^{(1)}=0$ and $X_{\lambda_n-L_1+1}^{(1)}=\cdots=X_{\lambda_n}^{(1)}=1$, then $\theta_{\lambda_n}^{(1)}=\frac{\sum_{i=L_1}^{\infty}(i-L_1)p_{0,i}^{(1)}}{\sum_{i=L_1}^{\infty}p_{0,i}^{(1)}}$.

Choose $k_2$ so large and $\delta_2$ so small such that

$$\frac{k_2\delta_2}{\beta_1+\delta_2}-\frac{\delta_2\sum_{i=L_1}^{\infty}ip_{0,i}^{(1)}}{\beta_1(\beta_1+\delta_2)}>2,$$

$$k_2\delta_2<\frac{1}{100^2}$$

and

$$\sum_{(x_0,\dots,x_{N_2})\in\{0,1\}^{N_2}} |P(X_0^{(1)}=x_0,\dots,X_{N_2}^{(1)}=x_{N_2}|X_0^{(1)}=0)$$

$$-P(X_0^{(2)}=x_0,\dots,X_{N_2}^{(2)}=x_{N_2}|X_0^{(2)}=0)|\le \frac{1}{1000^2}.$$

In this way, for the $\mathcal{M}^{(2)}$ process for some $L_0<\lambda_{n_1}<N_1$ and for some other $L_1<\lambda_{n_2}<N_2$ the estimate $h_n(X_0^{(2)},\dots,X_{\lambda_n}^{(2)})$ will be smaller than the target by at least 1 with probability $1-\frac{2}{1000}-\frac{2}{1000^2}$. (Note that $\sum_{i=0}^\infty ip_{0,i}^{(1)}\ge\sum_{i=0}^\infty ip_{0,i}^{(0)}$.)

Inductively, assume at stage $j$ we have a Markov chain $\mathcal{M}^{(j)}$ which satisfies the conditions $(C_j)$:

There are integers $L_0<N_1<L_1<\cdots<N_j$ such that

$$P\Bigg(\bigcup_{n_1=1}^\infty,\dots,\bigcup_{n_j=1}^\infty\bigcap_{i=1}^j\Bigg\{L_{i-1}<\lambda_{n_i}<N_i, X^{(j)}_{\lambda_{n_i}-L_{i-1}}=0,$$

$$X^{(j)}_{\lambda_{n_i}-L_{i-1}+1}=\cdots=X^{(j)}_{\lambda_{n_i}}=1,$$

$$h_{n_i}(X_0^{(j)},\dots,X^{(j)}_{\lambda_{n_i}})<\frac{\sum_{h=L_{i-1}}^\infty(h-L_{i-1})p_{0,h}^{(j)}}{\sum_{h=L_{i-1}}^\infty p_{0,h}^{(j)}}-1\Bigg\}\Bigg|$$

$$X_0^{(j)}=0\Bigg)$$

$$>1-\sum_{i=1}^j\frac{2}{1000^i}, \tag{15}$$

$$\sum_{(x_0,\dots,x_{N_j})\in\{0,1\}^{N_j}} |P(X_0^{(j-1)}=x_0,\dots,X_{N_j}^{(j-1)}=x_{N_j}|X_0^{(j-1)}=0)$$

$$-P(X_0^{(j)}=x_0,\dots,X_{N_j}^{(j)}=x_{N_j}|X_0^{(j)}=0)|\le\frac{1}{1000^j} \tag{16}$$

and

$$\sum_{h=0}^\infty hp_{0,h}^{(j)}\le 1+\sum_{h=1}^j\frac{1}{100^h}. \tag{17}$$

Now we will define $\mathcal{M}^{(j+1)}$. For an arbitrary $N_j<L_j$ let $\sum_{i\ge L_j}p_{0,i}^{(j)}=\beta$. For some $\delta<0.25p_{0,0}^{(j)}$ and $k$ which will be specified later, let $p_{0,0}^{(j)}=p_{0,0}^{(j)}-\delta$ and

$p_{0,k}^{(j)} = p_{0,k}^{(j)} + \delta$. Now the change in

$$\frac{\sum_{i=L_j}^{\infty}(i-L_j)p_{0,i}^{(j+1)}}{\sum_{i=L_j}^{\infty}p_{0,i}^{(j+1)}} - \frac{\sum_{i=L_j}^{\infty}(i-L_j)p_{0,i}^{(j)}}{\sum_{i=L_j}^{\infty}p_{0,i}^{(j)}} = \frac{k\delta}{\beta+\delta} - \frac{\sum_{i=L_j}^{\infty} ip_{0,i}^{(j)}\delta}{\beta(\beta+\delta)}$$

and in

$$\sum_{i=1}^{\infty} ip_{0,i}^{(j+1)} - \sum_{i=1}^{\infty} ip_{0,i}^{(j)} = k\delta.$$

Now choose $L_j$ such that $3\beta < 100^{-(j+1)}$. Choose $N_{j+1}$ so big such that

$$P\Bigg(\bigcup_{n=1}^{\infty}\Bigg\{L_j < \lambda_n < N_{j+1}, X_{\lambda_n - L_j}^{(j)} = 0, X_{\lambda_n - L_j+1}^{(j)} = \cdots = X_{\lambda_n}^{(j)} = 1,$$
$$\left|h_n(X_0^{(j)},\ldots,X_{\lambda_n}^{(j)}) - \frac{\sum_{i=L_j}^{\infty}(i-L_j)p_{0,i}^{(j)}}{\sum_{i=L_j}^{\infty}p_{0,i}^{(j)}}\right| < \frac{1}{100}\Bigg\} \,\Big|\, X_0^{(j)} = 0\Bigg)$$
$$> 1 - \left(\frac{1}{1000}\right)^j.$$

Note that if $k > K = \max_{0\le i<j} \frac{\sum_{h=L_i}^{\infty} hp_{0,h}^{(j)}}{\sum_{h=L_i}^{\infty} p_{0,h}^{(j)}}$, then for all $0 \le i < j$,

$$\frac{\sum_{h=L_i}^{\infty} hp_{0,h}^{(j+1)}}{\sum_{h=L_i}^{\infty} p_{0,h}^{(j+1)}} \ge \frac{\sum_{h=L_i}^{\infty} hp_{0,h}^{(j)}}{\sum_{h=L_i}^{\infty} p_{0,h}^{(j)}}. \tag{18}$$

Choose $k > K$ so large and $\delta$ so small such that

$$\frac{k\delta}{\beta+\delta} - \frac{\delta\sum_{i=L_j}^{\infty} ip_{0,i}^{(j)}}{\beta(\beta+\delta)} > 2,$$
$$k\delta < \frac{1}{100^{j+1}}$$

and

$$\sum_{(x_0,\ldots,x_{N_{j+1}})\in\{0,1\}^{N_{j+1}}} |P(X_0^{(j)} = x_0,\ldots,X_{N_{j+1}}^{(j)} = x_{N_{j+1}}|X_0^{(j)} = 0)$$
$$- P(X_0^{(j+1)} = x_0,\ldots,X_{N_{j+1}}^{(j+1)} = x_{N_{j+1}}|X_0^{(j+1)} = 0)|$$
$$\le \frac{1}{1000^{j+1}}.$$

The resulting Markov chain $\mathcal{M}^{(j+1)}$ is irreducible and positive recurrent and so it yields a stationary and ergodic distribution and the inductive assumption holds for $j+1$.

Define $p_{i,j}^{(\infty)} = \lim_{n\to\infty} p_{i,j}^{(n)}$. The resulting Markov chain $\mathcal{M}^{(\infty)}$ is clearly irreducible and positive recurrent and so it yields a stationary and ergodic distribution. Let $X_n = 0$ if $M_n^{(\infty)} = 0$ and 1 otherwise. Clearly, by the induction and (18),

$$P\Bigg(\bigcup_{n_1=1}^{\infty},\dots,\bigcup_{n_j=1}^{\infty}\bigcap_{i=1}^{j}\Big\{L_{i-1}<\lambda_{n_i}<N_i, X_{\lambda_{n_i}-L_i}=0,$$
$$X_{\lambda_{n_i}-L_i+1}=\dots=X_{\lambda_{n_i}}=1,$$
$$h_{n_i}(X_0,\dots,X_{\lambda_{n_i}}) < \frac{\sum_{h=L_i}^{\infty}(h-L_i)p_{0,h}^{(\infty)}}{\sum_{h=L_i}^{\infty} p_{0,h}^{(\infty)}} - 1\Big\}\Big| X_0 = 0\Bigg)$$
$$> 1 - \sum_{i=1}^{\infty}\frac{2}{1000^i}.$$

Since the set (event) is decreasing in $j$ so

$$P\Bigg(\bigcup_{l=1}^{\infty}\bigcup_{n_l=1}^{\infty}\bigcap_{i=1}^{\infty}\Big\{L_{i-1}<\lambda_{n_l}<N_i, X_{\lambda_{n_l}-L_i}=0,$$
$$X_{\lambda_{n_l}-L_i+1}=\dots=X_{\lambda_{n_l}}=1,$$
$$h_{n_l}(X_0,\dots,X_{\lambda_{n_l}}) < \frac{\sum_{h=L_i}^{\infty}(h-L_i)p_{0,h}^{(\infty)}}{\sum_{h=L_i}^{\infty} p_{0,h}^{(\infty)}} - 1\Big\}\ \Big|\ X_0 = 0\Bigg)$$
$$\geq 1 - \sum_{i=1}^{\infty}\frac{2}{1000^i}$$

and

$$\sum_{h=0}^{\infty} h p_{0,h}^{(\infty)} \leq 1 + \sum_{h=1}^{\infty}\frac{1}{100^h}.$$

The proof of Theorem 4 is complete.

**7. Proof of Theorem 5.** Consider the largest $L$ such that

$$P(\tau(X_{-\infty}^0) < L) \leq 1 - \frac{\varepsilon}{2}.$$

Applying the ergodic theorem we get that almost surely,

$$\liminf_{n\to\infty}\frac{n}{\lambda_n^{(\varepsilon)}} \geq P(\tau(X_{-\infty}^0)\leq L) \geq 1 - \frac{\varepsilon}{2} > 1-\varepsilon.$$

It is also clear that

$$\limsup_{n\to\infty}\frac{n}{\lambda_n^{(\varepsilon)}} \leq P(\tau(X_{-\infty}^0)\leq L),$$

and so eventually we are predicting for finitely many blocks of 1's and by ergodicity the consistency of the estimator $h_n^{(\varepsilon)}(X_0,\dots,X_{\lambda_n^{(\varepsilon)}})$ is also established. Since $p_l^{(\varepsilon)}(X_0^{\lambda_n^{(\varepsilon)}})$ is a probability distribution, now by ergodicity its consistency in total variation follows immediately for the same reason and the proof of Theorem 5 is complete.

**Acknowledgment.** We thank the referees of an earlier version for several helpful remarks including the reference to Petrov's book.

MTA-BME Stochastics Research Group
Institute of Mathematics
Egry József utca 1
Building H
Budapest 1111
Hungary
E-mail: morvai@math.bme.hu

Hebrew University of Jerusalem
Institute of Mathematics
Jerusalem 91904
Israel
E-mail: weiss@math.huji.ac.il